\documentclass{article}
\usepackage{eurosym}
\usepackage{amsfonts}
\usepackage{amssymb}
\usepackage{epsfig}
\usepackage{graphicx}
\usepackage{epstopdf}
\usepackage{wrapfig}
\usepackage[font=small,labelfont=bf]{caption} 

\def\iint{{\int\!\!\!\!\int\!\!}}

\newtheorem{theorem}{Theorem}[section]
\newtheorem{proposition}{Proposition}[section]
\newtheorem{lemma}{Lemma}[section]

\newtheorem{remark}{Remark}[section]

\def\theequation{\thesection.\arabic{equation}}
\begin{document}

\title{On the exact number of monotone solutions of a simplified Budyko
climate model and their different stability. }
\author{Sabri BENSID and Jes\'{u}s Ildefonso D\'{I}AZ}
\maketitle

\begin{abstract}
We consider a simplified version of the Budyko diffusive energy balance
climate model. We obtain the exact number of monotone stationary solutions
of the associated discontinuous nonlinear elliptic with absorption. We show
that the bifurcation curve, in terms of the solar constant parameter, is
S-shaped. We prove the instability of the decreasing part and the stability
of the increasing part of the bifurcation curve. In terms of the Budyko
climate problem the above results lead to an important qualitative
information which is far to be evident and which seems to be new in the
mathematical literature on climate models. We prove that if the solar
constant is represented by $\lambda \in (\lambda _{1},\lambda _{2}),$ for
suitable $\lambda _{1}<\lambda _{2},$ then there are exactly two stationary
solutions giving rise to a free boundary (i.e. generating two symmetric
polar ice caps: North and South ones) and a third solution corresponding to
a totally ice covered Earth. Moreover, we prove that the solution with
smaller polar ice caps is stable and the one with bigger ice caps is
unstable.
\end{abstract}

\noindent {\textbf{AMS Classification:}} {{35B35}, 35J61, 35P30, 35K58, }%
86A10\newline
\newline
\noindent {\textbf{Keywords:}} Nonlinear eigenvalue problem, {discontinuous
nonlinearity, S-shaped bifurcation curve, stability, free boundary, energy
balance, Budyko climate model.}\hfill \break

\section{Introduction}

The main goal of this paper is double: in a first step we study the exact
number of solutions, depending on the parameter $\lambda $, of the
discontinuous eigenvalue type problem
\[
P(\lambda ,f)\left\{
\begin{array}{ll}
-u_{xx}(x)+\omega ^{2}u(x)=\lambda f(u(x)) & \quad \mbox{
in }\ x\in (0,1), \\[0.3cm]
u_{x}(0)=0, u(1)=0, & \quad%
\end{array}%
\right.
\]%
where $\omega ^{2}$ is a given parameter and $f(u)$ is the discontinuous
function given by
\begin{equation}
f(v)=f_{0}+(1-f_{0})H(v-\mu ),  \label{f(u)}
\end{equation}%
for some $\mu >0,$ under the key assumption
\begin{equation}
f_{0}\in (0,1),  \label{f_0}
\end{equation}%
with $H(s)$ the Heaviside discontinuous function
\[
H(s)=0\quad \hbox{for}\quad s<0,\quad \quad H(s)=1\quad \hbox{for}\quad
s\geq 0.
\]

Our second and main motivation is to study the different stability of the
monotone solutions of $P(\lambda ,f)$ as stationary solutions associated to
the parabolic problem
\[
PP^{\ast }(\lambda ,\beta ,u_{0})\left\{
\begin{array}{ll}
u_{t}-u_{xx}+\omega ^{2}u\in \lambda \beta (u) & \quad x\in (0,1),t>0, \\%
[0.3cm]
u_{x}(0,t)=0,\quad u(1,t)=0 & \quad t>0, \\[0.3cm]
u(x,0)=u_{0}(x), & \quad x\in (0,1).%
\end{array}%
\right.
\]%
Here $\beta $ is the maximal monotone graph of $\mathbb{R}^{2}$ given by
\[
\left\{
\begin{array}{ll}
\beta (r)=f(r) & \hbox{if }r\neq \mu , \\[0.3cm]
\beta (\mu )=[f_{0},1]. &
\end{array}%
\right.
\]%
Notice that any solution of $P(\lambda ,f)$ is also a weak solution of the
multivalued problem
\[
P^{\ast }(\lambda ,\beta )\left\{
\begin{array}{ll}
-u^{\prime \prime }(x)+\omega ^{2}u\in \lambda \beta (u(x)) & a.e.\hbox{ }%
x\in (0,1), \\[0.3cm]
u^{\prime }(0)=0 \hbox{, }u(1)=0. & \quad%
\end{array}%
\right.
\]

The present paper can be considered as a natural continuation of the
previous paper by the authors \cite{BensidDiaz} in which the same program of
research was devoted to the special case without any absorption term $\omega
=0$. As we shall see, several of the sharp methods of proof used in \cite%
{BensidDiaz} do not admit any easy adaptation to the apparently minor change
introduced when assuming $\omega ^{2}>0.$

As mentioned in \cite{BensidDiaz}, Problem $P(\lambda ,f)$ can be considered
as a simplified version of some more general formulations arising in several
different contexts: chemical reactors and porous media combustion, steady
vortex rings in an ideal fluid, plasma studies, the primitive equations of
the atmosphere in presence of vapor saturation, etc. We send the reader to
the references collected in \cite{BensidDiaz}. Nevertheless, our special
motivation to study the case $\omega ^{2}>0$ was the consideration of
problem $P(\lambda ,f)$ as a simplified version of the so called Budyko
diffusive energy balance models arising in climatology (see, e.g. \cite%
{North}, \cite{Hetzer Houston}, \cite{Diaz Escorial}, \cite{Stakgold}, \cite%
{Diaz Tello Collecatanea} and a stochastic version in \cite{Diaz Langa
Valero}). Although these models must be formulated on a Riemannian manifold
without boundary representing the Earth atmosphere \cite{Diaz Tello
Collecatanea}, the so called 1d-model corresponds to the case in which the
surface temperature is assumed to depend only on the latitude component. The
model considered in our previous paper \cite{BensidDiaz} neglected the
important \textit{feedback term} arising when modeling the emitted
terrestrial energy flux ( represented here, in a simplified form, by the
term $\omega ^{2}u-f_{0}$). In this way, we lead to a formulation similar to
$P(\lambda ,f)$ in which the spatial domain $(0,1)$ must be associated to a
semisphere, the discontinuous function represents the co-albedo (with a
discontinuity which is associated to the radical change of the co-albedo
when the temperature is crossing $-10$ centigrade degrees: here represented
by the value $u=\mu $), the parameter $\lambda $ the so-called solar
constant, the boundary condition $u^{\prime }(0)=0$ formulates the
simplified assumption of symmetry between both semispheres and the condition
$u(1)=0$ represents the renormalized temperature at the North pole (i.e. we
are assuming that $u\geq 0$ in the rest of the hemisphere, and thus we must
assume that $\mu >0$ although it represents $-10$ centigrade degrees).

We point out that the case in which the absorption term $\omega ^{2}u$ plays
also an important role in the class of eigenvalue type problems $P(\lambda
,f)$ associated to a discontinuous nonlinearity corresponds to the modelling
considered by McKean \cite{MKean} of the initial value problem for the
FitzHugh--Nagumo equations which were introduced as a model for the
conduction of electrical impulses in the nerve axon (see, e.g., Terman \cite%
{Termam}).

We also recall that the very sharp bifurcation and stability results
obtained trough the famous Crandall-Rabinowitz paper \cite%
{Crandall-Rabinowitz} requires in a fundamental way the differentiability of
the nonlinear term. That was used in the very nice paper \cite{Hetzer
Houston} to study the Sellers diffusive energy balance climate model in
which $\beta (u)$ is assumed to be at least a Lipschitz continuous function.

Results on the asymptotic behavior, when $t\rightarrow +\infty $, for the
evolution energy balance model were obtained in \cite{Diaz-Her-Tello} where
it was also proved the general multiplicity of stationary solutions
according the value of $\lambda $ (see also \cite{Valero} and \cite{Gorban}
for other related results)$.$ A sharper bifurcation diagram, as a S-shaped
curve was rigorously obtained in \cite{ADT}. Nevertheless the method of
proof in \cite{ADT} uses the information obtained trough suitable
zero-dimensional energy balance models and thus there is lacking of a more
detailed information about the associated free boundaries generated by the
solutions (given as the spatial points where $T=-10$). This kind of sharper
information will be obtained here since no zero-dimensional energy balance
model will be used in our proofs but only a direct analysis of the 1d-model.

One of the main difficulties to adapt the tools used in our previous paper
\cite{BensidDiaz} to the case in which $\omega ^{2}$ is not zero is the fact
that the solutions of the ODE $-u^{\prime \prime }+\omega ^{2}u=M$ may
oscillate (in contrast with the case $\omega ^{2}=0$). As a matter of fact,
there are several results in the literature indicating that the Budyko
diffusive energy balance climate model admits an infinity of stationary
solutions. That was shown in the \cite{Schmidt Bertina} and \cite{Hetzer
infinite} by considering a non-autonomous term $\lambda f(x,u)$ and in \cite%
{Diaz Tello infinite} for the mere autonomous case. The main goal of this
paper concerns the study of non-oscillating solutions of the autonomous
framework (which in the title is referred as \textquotedblleft monotone
solutions\textquotedblright , as we shall explain below).

In order to state our results we start by defining two crucial values of the
parameter $\lambda $:
\[
\lambda _{1}:=\mu \frac{\omega ^{2}\cosh \omega }{\sinh (\omega \kappa
)\sinh (\omega -\omega \kappa )-f_{0}\cosh (\omega \kappa )+f_{0}\cosh
(\omega \kappa )\cosh (\omega -\omega \kappa )},
\]%
and%
\[
\lambda _{2}:=\mu \frac{\omega ^{2}\cosh \omega }{f_{0}(\cosh \omega -1)},
\]%
where $\kappa =\kappa (f_{0})\in (0,1)$ will be given later (see formula (%
\ref{kappa}) below and take $\kappa =r^{\ast }$)$.$ We shall use the
notation $\Vert u\Vert _{\infty }=\displaystyle\max_{x\in \lbrack
0,1]}|u(x)| $.

By a solution $u_{\lambda ,\mu }$ of problem $P(\lambda ,f)$ we mean a
function $u\in C^{2}((0,1)\setminus \{x_{\lambda,\mu}\})\cap C^{1}([0,1))$,
for some $x_{\lambda,\mu }\in \lbrack 0,1)$ where $u(x_{\lambda,\mu })=\mu $
(called as \textit{the free boundary associated to} $u$) and with $u\geq 0$,
$u\neq 0$, such that $-u^{\prime \prime }(x)+\omega ^{2}u(x)=\lambda
f(u(x)), $ for any $x\in (0,1)-\{x_{\lambda,\mu }\},$ and \ $u^{\prime
}(0)=0,$ $u(1)=0$. \ As mentioned before, our main interest concerns \textit{%
monotone solutions} $u_{\lambda ,\mu }$ of problem $P(\lambda ,f)$ \ (i.e.
such that, in addition, $u^{\prime }(x)\leq 0$ for any $x\in (0,1)$). \ Our
first result shows the exact multiplicity of \textit{monotone solutions} $%
u_{\lambda ,\mu }$ of problem $P(\lambda ,f)$ for different values of $%
\lambda $.

\begin{theorem}
i) If $\lambda <\lambda _{2},$ then there exists a unique solution $%
u_{\lambda ,\mu }^{\ast }$ without free boundary of $P(\lambda ,f)$.
Moreover $u_{\lambda ,\mu }^{\ast }$ is monotone and
\begin{equation}
\Vert u_{\lambda ,\mu }^{\ast }\Vert _{\infty }=u_{\lambda ,\mu }^{\ast
}(0)=-\frac{\lambda f_{0}}{\omega^2\cosh \omega }+\frac{\lambda f_{0}}{%
\omega^2}<\mu ,  \label{u star maximum}
\end{equation}%
i.e. the line $(\lambda ,\gamma ^{\ast }(\lambda ))$
\[
\gamma ^{\ast }(\lambda ):=\frac{\lambda f_{0}(\cosh \omega -1)}{%
\omega^2\cosh \omega },\hbox{ if }\lambda \in (0,\lambda _{2}),
\]%
defines an increasing part of the $\lambda -$bifurcation diagram.

\noindent ii) If $\lambda =\lambda _{1}$ then there exists a unique monotone
solution $u_{\lambda _{1}}$ of $P(\lambda ,f)$ giving rise to a free
boundary. Moreover $u_{\lambda _{1}}$ is strictly concave and $u_{\lambda
_{1}}(0)=\mu .$

\noindent iii) If $\lambda \in (\lambda _{1},\lambda _{2}]$ then there
exists $\underline{u}_{\lambda ,\mu }$ monotone solution of $P(\lambda ,f)$
with a free boundary $\underline{x}_{\lambda,\mu }\in (0,r^{\ast }),$ where $%
r^{\ast }\in (0,1)$ is such that $u_{\lambda ,\mu }(r^{\ast })=\mu .$
Moreover,
\[
\Vert \underline{u}_{\lambda ,\mu }\Vert _{\infty }=\underline{u}_{\lambda
,\mu }(0)=\frac{(\mu\omega^2 -\lambda) }{\omega^2\cosh (\omega \underline{x}%
_{\lambda,\mu })}+\frac{\lambda}{\omega^2} :=\underline{\gamma }(\lambda ),
\]

\noindent iv) If $\lambda \in (\lambda _{1},+\infty )$ then there exists $%
\overline{u}_{\lambda ,\mu }$ monotone solution of $P(\lambda ,f)$ such that
$\mu <\Vert \overline{u}_{\lambda ,\mu }\Vert _{\infty }$ and $\Vert
\underline{u}_{\lambda ,\mu }\Vert _{\infty }<\Vert \overline{u}_{\lambda
,\mu }\Vert _{\infty }$ if $\lambda \in (\lambda _{1},\lambda _{2}].$
Moreover its free boundary is given by $\overline{x}_{\lambda,\mu}\in
(r^{\ast },1)$, with $r^{\ast }\in (0,1)$ given as in iii), and
\[
\Vert \overline{u}_{\lambda ,\mu }\Vert _{\infty }=\overline{u}_{\lambda
,\mu }(0)=\frac{(\mu \omega^2-\lambda) }{\omega^2\cosh (\omega \overline{x}%
_{\lambda,\mu })}+\frac{\lambda}{\omega^2} :=\overline{\gamma }(\lambda ).
\]

v) The $\lambda -$bifurcation curve is S-shaped, i.e. it is a continuous
curve of $\lambda $ such that $\gamma ^{\ast }(\lambda )$, $\overline{\gamma
}(\lambda )$ (respectively $\underline{\gamma }(\lambda )$) are increasing
(respectively decreasing) functions of $\lambda .$
\end{theorem}

\begin{figure}[h]
\label{Bifucation S-shaped curve}
\begin{center}
\includegraphics[scale=0.8]{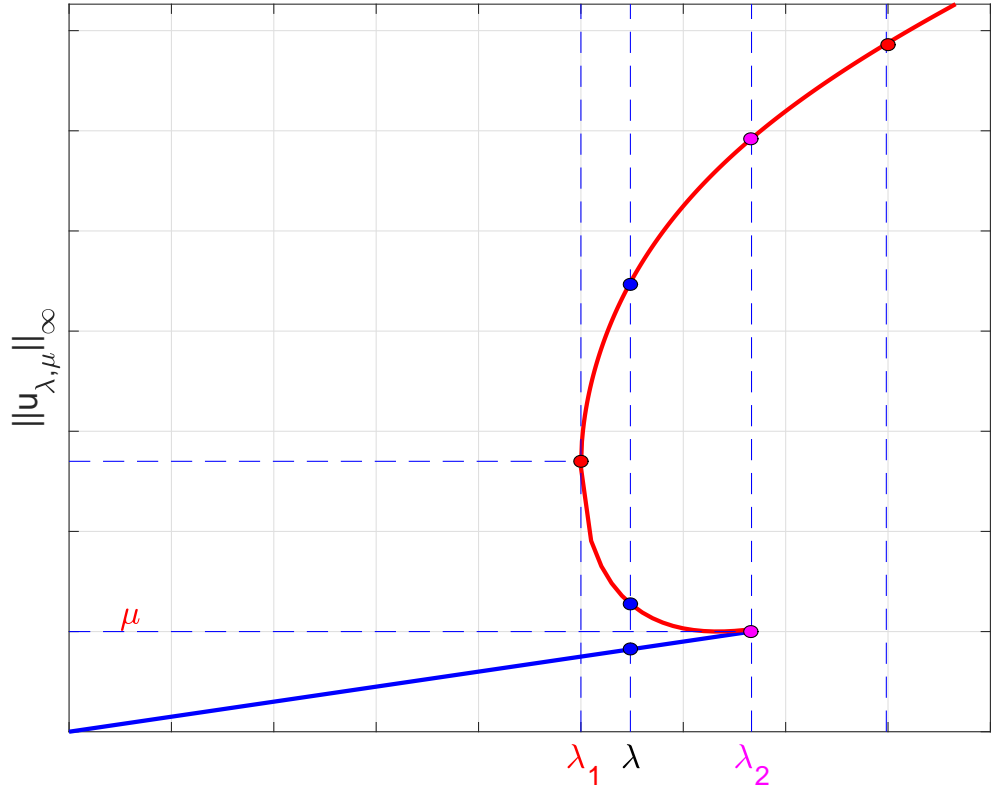}
\caption{Bifucation S-shaped curve.}
\end{center}
\end{figure}
In order to state our results concerning the stability of the above monotone
solutions of $P(\lambda ,f)$\ as stationary solutions of the associated
parabolic problem $PP^{\ast }(\lambda ,\beta ,u_{0})$, we start by denoting
by $u(t,x:u_{0})$ to the solution of $PP^{\ast }(\lambda ,\beta ,u_{0}).$ It
is a routine matter to check that all the existence and uniqueness results
presented in \cite{BensidDiaz} for problem $PP^{\ast }(\lambda ,\beta
,u_{0}) $ for the case $\omega =0$ extends without difficulty to the case $%
\omega ^{2}>0$ (many references on the previous literature on problems
related to $PP^{\ast }(\lambda ,\beta ,u_{0})$ were presented in \cite%
{BensidDiaz}; see also \cite{Deguchi1}). We recall that, due to the
discontinuity of $f(u)$, the uniqueness of solution of $PP^{\ast }(\lambda
,\beta ,u_{0})$\ (and the comparison principle) requires to work in the
class of \textquotedblleft non degenerate solutions": \ i.e. solutions $%
u(t,x:u_{0})$ such that
\[
meas\{x\in (0,1),\mid u(t,x:u_{0})-\mu \mid \leq \theta \}\leq C\theta ,
\]%
for any $\theta \in (0,\theta _{0})$ and for any $t>0$, for some $C>0$ and $%
\theta _{0}>0$. Here $meas(.)$ denotes the Lebesgue measure.

The following theorem shows that if $\lambda \in (\lambda _{1},+\infty ),$
then $\overline{u}_{\lambda ,\mu }$ is stable in $L^{\infty }(0,1).$

\begin{theorem}
Let $u_{0}\in L^{\infty }(0,1)$ with $u_{0}\geq 0$ a.e in $(0,1).$ Let $%
\overline{u}_{\lambda ,\mu }$ be the monotone solution of $P(\lambda ,f)$
given in $iv).$\newline
Assume that $\Vert u_{0}-\overline{u}_{\lambda ,\mu }\Vert _{L^{\infty }}$
is sufficiently small and let $u(t,x:u_{0})$ be the non-degenerate solution
of $PP^{\ast }(\lambda ,\beta ,u_{0})$. Then
\[
\Vert u(x,t:u_{0})-\overline{u}_{\lambda ,\mu }(x)\Vert _{L^{\infty
}(0,1)}<\varepsilon \quad \quad \hbox{\ for any }t\geq 0,
\]%
for some positive constant $\varepsilon .$
\end{theorem}

\bigskip

The proof of the instability of the stationary monotone solutions $%
\underline{u}_{\lambda ,\mu }$ of the decreasing part of the bifurcation
curve $(\lambda ,\underline{\gamma }(\lambda ))$ can be obtained by
different tools according the values of the parameter $\omega .$ A first
possibility (as in \cite{BensidDiaz}) is to study the sign of the first
eigenvalue of the problem associated to the linearized equation
\[
P_{\eta }(x_{\lambda ,\mu }:f_{0},\lambda )\left\{
\begin{array}{ll}
-U^{\prime \prime }(x)-\lambda (1-f_{0})\delta _{\{x_{\lambda ,\mu
}\}}U(x)+\omega ^{2}U(x)=\eta U(x), & \quad x\in (0,1) \\[0.3cm]
U^{\prime }(0)=0,\quad U(1)=0, & \quad%
\end{array}%
\right.
\]%
where $\delta _{\{x_{\lambda ,\mu }\}}$ is the Dirac delta distribution at
the free boundary point $x_{\lambda ,\mu }.$ When $0<\omega <1$ we can adapt
the study made in \cite{BensidDiaz} to prove that the principal eigenvalue $%
\eta _{1}$ of problem $P_{\eta }(x_{\lambda ,\mu }:f_{0},\lambda )$ is
negative.

\begin{theorem}
If $0<\omega <1$ and $\lambda \in (\lambda _{1},\lambda _{2}],$ then the
stationary monotone solution $\underline{u}_{\lambda ,\mu }$ is unstable in $%
L^{\infty }(0,1).$
\end{theorem}

Unfortunately, when $\omega \geq 1$ the above method of proof is not
applicable and other arguments are needed. The following results use the
sharp description of the equilibria given in Theorem 1.1 jointly to some
monotonicity arguments$.$

\begin{theorem}
Let $\lambda \in (\lambda _{1},\lambda _{2}].$ For any $\varepsilon >0,$
there exists $\underline{u}_{0},\overline{u}_{0}\in L^{\infty }(0,1)$, non
degenerate functions such that\newline
\[
\begin{array}{ll}
\underline{u}_{\lambda ,\mu }-\underline{u}_{0}\leq 0, & \underline{u}%
_{\lambda ,\mu }-\overline{u}_{0}\geq 0,%
\end{array}%
\]%
and
\begin{equation}
\begin{array}{ll}
||\underline{u}_{\lambda ,\mu }-\underline{u}_{0}||_{L^{\infty
}}<\varepsilon , & ||\underline{u}_{\lambda ,\mu }-\overline{u}%
_{0}||_{L^{\infty }}<\varepsilon .%
\end{array}
\label{epsilon near}
\end{equation}%
Moreover, if $u(t,x:\overline{u}_{0})$ (respectively $u(t,x:\underline{u}%
_{0})$) is the unique non degenerate solution of the corresponding problem $%
PP^{\ast }(\lambda ,\beta ,u_{0}),$ then
\begin{equation}
\frac{\partial u}{\partial t}(t,.:\overline{u}_{0})\geq 0\quad \hbox{a.e}%
\quad t>0
\end{equation}%
\begin{equation}
u(t,.:\overline{u}_{0})\leq \overline{u}_{\lambda ,\mu }(.)\quad \forall
t\in \lbrack 0,+\infty )
\end{equation}%
\begin{equation}
u(t,.:\overline{u}_{0})\rightarrow \overline{u}_{\lambda ,\mu }(.)\quad %
\hbox{in}\quad H^{1}(0,1)\quad \hbox{when}\quad t\rightarrow +\infty ,
\end{equation}%
respectively
\begin{equation}
\frac{\partial u}{\partial t}(t,.:\underline{u}_{0})\leq 0\quad \hbox{a.e}%
\quad t>0
\end{equation}%
\begin{equation}
u(t,.:\underline{u}_{0})\geq u_{\lambda ,\mu }^{\ast }(.)\quad \forall t\in
\lbrack 0,+\infty )  \label{convergence u star 0}
\end{equation}%
\begin{equation}
u(t,.:\underline{u}_{0})\rightarrow u_{\lambda ,\mu }^{\ast }(.)\quad %
\hbox{in}\quad H^{1}(0,1)\quad \hbox{when}\quad t\rightarrow +\infty .
\label{convergence u star}
\end{equation}%
In particular, the stationary monotone solution $\underline{u}_{\lambda ,\mu
}(.)$ is unstable in $H^{1}$ (and so unstable also in $L^{\infty }$).
\end{theorem}

In terms of the Budyko climate problem the above theorems lead to an
important qualitative information which is far to be evident and which seems
to be new in the mathematical literature on climate models. Let us extend,
just by symmetry, the solutions of $P(\lambda ,f)$ and $PP^{\ast }(\lambda
,\beta ,u_{0})$ to the whole spatial interval $(-1,1)$ (this corresponds to
consider the atmospheric surface temperature on the whole sphere instead on
the south hemisphere). Assume the solar constant $\lambda \in (\lambda
_{1},\lambda _{2}).$ Then there are exactly two stationary solutions giving
rise to a free boundary (i.e. generating two symmetric polar ice caps: North
and South ones) and a third solution corresponding to a totally ice covered
Earth. Moreover, among the two more realistic solutions (presenting two
symmetric polar ice caps) the solution with smaller polar ice caps is stable
and the one with bigger ice caps is unstable. \newline
\newline

\begin{center}
\includegraphics[scale=0.35]{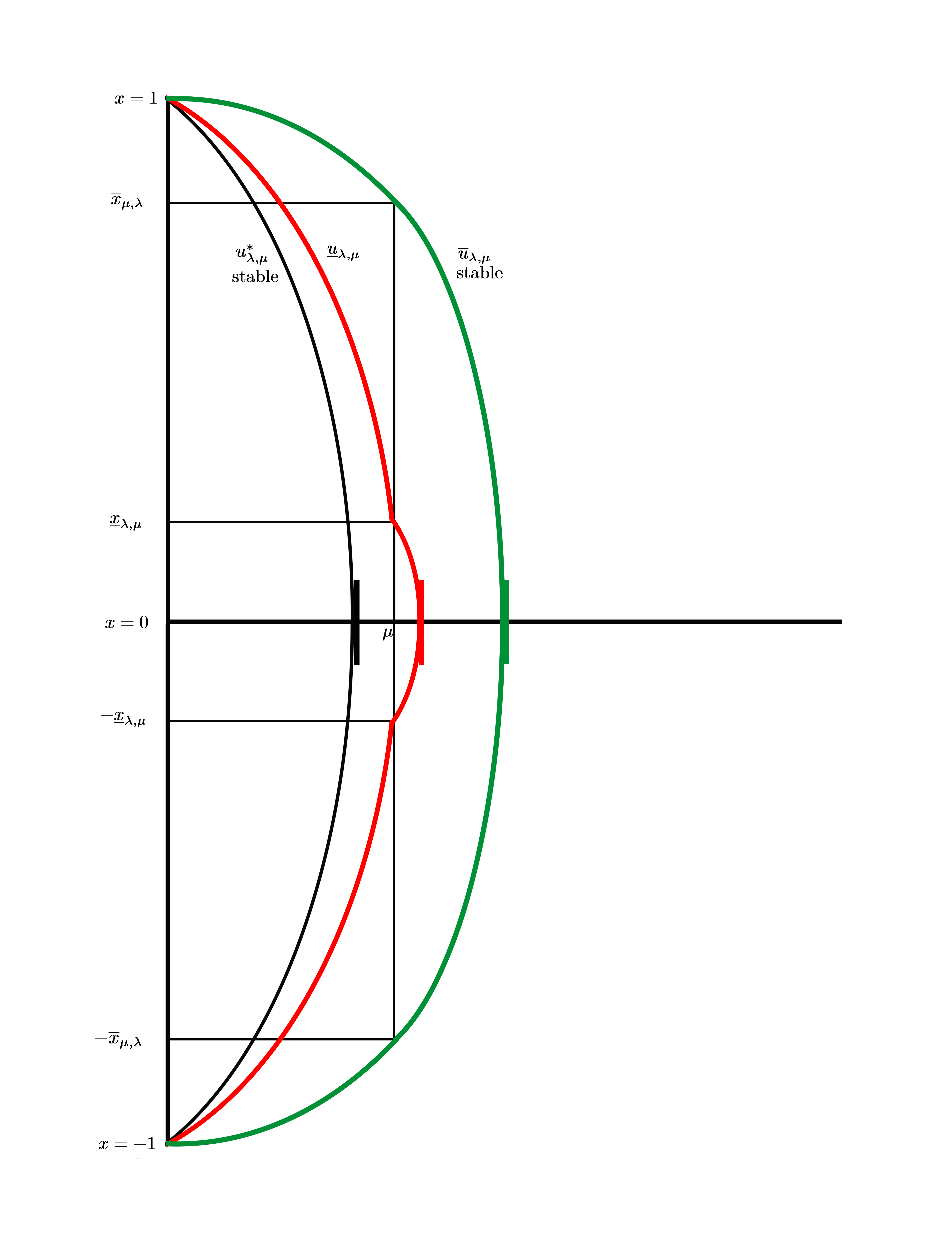}
\captionof{figure}{Qualitative representation of the three asthenosphere equilibria
temperature depending of the equilatitude parallel circles $x\in \lbrack
-1,1]$. }
\end{center}

Obviously, there are many aspects of the correct modeling of the energy
balance climate models which were not taken into account in the simpler
model analyzed in this paper. Nevertheless, it seems difficult to imagine
that the stability properties of the more complex states which correspond to
the monotone stationary solutions studied in our framework may be completely
different to what is suggested here thanks to the \textquotedblleft
simplicity\textquotedblright\ of the model.

\section{Proof of Theorem 1.1: the S-shaped bifurcation diagram}

To prove i) we must consider the case where $f(u(x))=f_{0}$ for any $x\in
\lbrack 0,1]$ (the case without free boundary)$.$ Hence, we have
\[
P(f_{0})\left\{
\begin{array}{ll}
-u^{\prime \prime }+\omega ^{2}u=\lambda f_{0} & \quad \mbox{
in }\ (0,1), \\[0.3cm]
u^{\prime }(0)=0,\hbox{ }u(1)=0, & \quad%
\end{array}%
\right.
\]%
An easy calculation (see, e.g. \cite{zill}) shows that the solution of of $%
P(f_{0})$ is given by
\[
u^{\ast }(x)=-\frac{\lambda f_{0}}{\omega^2\cosh \omega }\cosh (\omega x)+%
\frac{\lambda f_0}{\omega^2},\quad \hbox{for any }\quad x\in \lbrack 0,1].
\]%
Since $u^{\ast }(0)=-\frac{\lambda f_{0}}{\omega^2\cosh \omega }+\frac{%
\lambda f_{0}}{\omega^2}<\mu ,$ we have $u^{\ast }(0)<\mu $ if and only if
\[
\frac{\lambda }{\mu }<\frac{\omega^2\cosh \omega }{f_{0}(\cosh \omega -1)}.
\]%
Now, we shall search monotone solutions $u$ with a free boundary $x_{\mu
,\lambda }\in (0,1)$ and we consider the corresponding problems verified by $%
u$ on the different regions $(0,x_{\lambda,\mu })$ and $(x_{\lambda,\mu
},1). $ On $(0,x_{\lambda,\mu }),$ we get the following problem
\[
(P_{L})\left\{
\begin{array}{ll}
-u^{\prime \prime }+\omega ^{2}u=\lambda & \quad \mbox{
in }\ (0,x_{\lambda,\mu }), \\[0.3cm]
u^{\prime }(0)=0,\hbox{ }u(x_{\lambda,\mu })=\mu . & \quad%
\end{array}%
\right.
\]%
Then
\begin{equation}
u_{\lambda ,\mu }(x)=\frac{(\mu\omega^2 -\lambda )\cosh (\omega x)}{%
\omega^2\cosh (\omega x_{\lambda,\mu })}+\frac{\lambda}{\omega^2} .
\label{u  on the left}
\end{equation}%
On $(x_{\lambda,\mu },1),$ we get
\[
(P_{R})\left\{
\begin{array}{ll}
-u^{\prime \prime }+\omega^2 u=\lambda f_{0} & \quad \mbox{
in }\ (x_{\lambda,\mu },1), \\[0.3cm]
u(x_{\lambda,\mu })=\mu ,\hbox{ }u(1)=0. & \quad%
\end{array}%
\right.
\]%
So,
\begin{eqnarray}
u_{\lambda ,\mu }(x) &=&\left[ -\frac{\lambda f_{0}}{\omega^2\cosh \omega }%
-\left( \frac{\mu \cosh \omega -\frac{\lambda f_{0}}{\omega^2}\cosh \omega +%
\frac{\lambda f_{0}}{\omega^2}\cosh (\omega x_{\lambda,\mu })}{\sinh (\omega
x_{\lambda,\mu }-\omega )}\right) \frac{\sinh \omega }{\cosh \omega }\right]
\cosh (\omega x) \\
&&+\left[ \frac{\mu \cosh \omega -\frac{\lambda f_{0}}{\omega^2}\cosh \omega
+\frac{\lambda f_{0}}{\omega^2}\cosh (\omega x_{\lambda,\mu })}{\sinh
(\omega x_{\lambda,\mu }-\omega )}\right] \sinh (\omega x)+\frac{\lambda
f_{0}}{\omega^2}.
\end{eqnarray}

\noindent The transmission condition lead to the necessary condition
\begin{equation}
\frac{\lambda }{\mu }=\frac{\omega ^{2}\cosh \omega }{\sinh (\omega
x_{\lambda,\mu })\sinh (\omega -\omega x_{\lambda,\mu })-f_{0}\cosh (\omega
x_{\lambda,\mu })+f_{0}\cosh (\omega x_{\lambda,\mu })\cosh (\omega -\omega
x_{\lambda,\mu })}  \label{TC}
\end{equation}%
In order to study this condition, let us introduce the auxiliary function
\begin{equation}
g(r)=\frac{\omega ^{2}\cosh \omega }{\sinh (\omega r)\sinh (\omega -\omega
r)-f_{0}\cosh (\omega r)+f_{0}\cosh (\omega r)\cosh (\omega -\omega r)}.
\label{function g}
\end{equation}%
Hence,
\[
g^{\prime }(r)=0\hbox{ if }\frac{-\omega ^{2}\cosh \omega \left[ \sinh
(\omega -2\omega r)+f_{0}\sinh (2\omega r-\omega )-f_{0}\sinh (\omega r)%
\right] }{(\sinh (\omega r)\sinh (\omega -\omega r)-f_{0}\cosh (\omega
r)+f_{0}\cosh (\omega r)\cosh (\omega -\omega r))^{2}}=0.
\]%
This implies that
\[
\sinh (2\omega r-\omega )-f_{0}\sinh (2\omega r-\omega )+f_{0}\sinh (\omega
r)=0.
\]%
Thus,
\begin{equation}
g^{\prime }(r^{\ast })=0\hbox{ if and only if }r^{\ast }=\frac{%
m_{f_{0}}(r^{\ast })+1}{2},  \label{kappa}
\end{equation}%
where
\[
m_{f_{0}}(r):=\frac{1}{\omega }\sinh ^{-1}\left( \frac{f_{0}}{f_{0}-1}\sinh
(\omega r)\right) .
\]%
Clearly, there is a unique fixed point $r^{\ast }\in (0,1)$ and the function
$g$ has a minimum equal to $\lambda _{1}$ at $r^{\ast }.$ Moreover, $g$ is
monotone decreasing in $(0,r^{\ast })$ and monotone increasing in $(r^{\ast
},1).$\newline
Hence, when $\lambda =\lambda _{1},$ it follows that equation (\ref{TC}) has
one root $x_{\lambda_1,\mu}$ and we obtain the desired monotone solution
noted by $u_{\lambda _{1}}.$ \newline
When, $\lambda >\lambda _{1},$ the equation (\ref{TC}) has two roots $%
\overline{x}_{\lambda,\mu }$ and $\underline{x}_{\lambda,\mu }$ between $%
(0,1)$ different from $x_{\lambda_1,\mu}.$\newline
If we denote by $\underline{u}_{\lambda ,\mu }(x)$ and $\overline{u}%
_{\lambda ,\mu }(x)$ the functions satisfying (2.10) and (2.12) and with
free boundaries given respectively by $\underline{x}_{\lambda,\mu }$ and $%
\overline{x}_{\lambda,\mu},$ then we get the conclusions stated in $iii)$
and $iv).$

\noindent To prove part v) of Theorem 1.1 we introduce the auxiliary
function
\[
h_{\varepsilon ,\omega }(x)=\cosh \omega -\sinh (\omega x)\sinh (\omega
-\omega x)+\varepsilon \cosh (\omega x)-\varepsilon \cosh (\omega x)\cosh
(\omega -\omega x)
\]%
\ \noindent for $x\in (0,1)$ and $\varepsilon \in (0,1)$. It is easy to see
that for $\varepsilon \rightarrow 1,$ $h_{\varepsilon ,\omega }(x)>0$ as
well as also for $\varepsilon \rightarrow 0.$ Since $h_{\varepsilon ,\omega
} $ is nondecreasing with respect to $\varepsilon ,$ then we conclude that $%
h_{\varepsilon ,\omega }(x)>0$ for any $x\in (0,1)$ and $\varepsilon \in
(0,1)$. Hence, using the fact that $h_{\varepsilon ,\omega }(x)>0,$ we
obtain
\[
\cosh \omega >\sinh (\omega x)\sinh (\omega -\omega x)-\varepsilon \cosh
(\omega x)+\varepsilon \cosh (\omega x)\cosh (\omega -\omega x).
\]%
Thus,%
\[
\frac{\omega ^{2}\cosh \omega }{\sinh (\omega r)\sinh (\omega -\omega
r)-f_{0}\cosh (\omega r)+f_{0}\cosh (\omega r)\cosh (\omega -\omega r)}%
>\omega ^{2}
\]%
giving that $\lambda _{1}>\mu \omega ^{2}.$ Note also that the bifurcation
branch $\underline{\gamma }$ is given by
\[
\underline{\gamma }(\lambda )=\frac{(\mu \omega ^{2}-\lambda )}{\omega
^{2}\cosh (\omega \underline{x}_{\lambda ,\mu })}+\frac{\lambda }{\omega ^{2}%
},\quad \hbox{for}\quad \lambda \in (\lambda _{1},\lambda _{2}].
\]%
When $\lambda \rightarrow \lambda _{2},$ $x_{\lambda ,\mu }\rightarrow
x_{\lambda _{2},\mu }=0.$\newline
Hence,
\[
\underline{\gamma }(\lambda )=\frac{(\mu \omega ^{2}-\lambda _{2})}{\omega
^{2}\cosh (0)}+\frac{\lambda _{2}}{\omega ^{2}}=\mu .
\]%
When $\lambda \rightarrow \lambda _{1},$ $x_{\lambda ,\mu }\rightarrow
x_{\lambda _{1},\mu }=r^{\ast }.$ In this case,
\[
\underline{\gamma }(\lambda _{1})=\frac{(\mu \omega ^{2}-\lambda _{1})}{%
\omega ^{2}\cosh (\omega ^{2}r^{\ast })}+\frac{\lambda _{1}}{\omega ^{2}}=%
\frac{\mu \omega ^{2}}{\omega ^{2}\cosh (\omega r^{\ast })}+\frac{\lambda
_{1}}{\omega ^{2}}(\frac{\omega ^{2}\cosh (\omega r^{\ast })-1}{\omega
^{2}\cosh (\omega r^{\ast })}).
\]%
Using the fact that $\lambda _{1}>\mu \omega ^{2},$ we have
\[
\underline{\gamma }(\lambda _{1})>\frac{\mu \omega ^{2}}{\omega ^{2}\cosh
(\omega r^{\ast })}+\mu (\frac{\omega ^{2}\cosh (\omega r^{\ast })-1}{\omega
^{2}\cosh (\omega r^{\ast })})>\mu .
\]%
The same result holds for the part of the bifurcation curve $\overline{%
\gamma }(\lambda )$ given by
\[
\overline{\gamma }(\lambda )=\Vert \overline{u}_{\lambda ,\mu }\Vert
_{\infty }=\overline{u}_{\lambda ,\mu }(0)=\frac{\mu \omega ^{2}-\lambda }{%
\omega ^{2}\cosh (\omega \overline{x}_{\lambda ,\mu })}+\frac{\lambda }{%
\omega ^{2}},\quad \hbox{for}\quad \lambda \in (\lambda _{1},+\infty ).
\]%
Notice that when $\lambda \rightarrow +\infty ,$ $\overline{\gamma }%
\rightarrow +\infty .$ The proof of Theorem 1.1 is then completed.$%
_{\blacksquare }$

\section{Stability of the increasing part of the bifurcation curve.}

The proof of Theorem 1.2 will follow the same philosophy than the proof of
Theorem 1.3 of \cite{BensidDiaz}, but with some important modifications
based on the two following results:

\begin{proposition}
Let the assumptions of Theorem 1.2 hold. Then, for some positive parameters $%
\theta $ and $\delta $, there exist two continuous functions $\overline{\psi
}_{\theta ,\delta }(x:\lambda ,\mu ),\underline{\psi }_{\theta ,\delta
}(x:\lambda ,\mu )$, which depend continuously of the parameters $\theta $
and $\delta $, and $\varepsilon =\varepsilon (\theta ,\delta )>0$, such that
\begin{equation}
\overline{u}_{\lambda ,\mu }(x)-\varepsilon \leq \underline{\psi }_{\theta
,\delta }(x:\lambda ,\mu )<\overline{u}_{\lambda ,\mu }(x)<\overline{\psi }%
_{\theta ,\delta }(x:\lambda ,\mu )\leq \overline{u}_{\lambda ,\mu
}(x)+\varepsilon \hbox{, for any }x\in \lbrack 0,1].  \label{f2}
\end{equation}
In particular,
\begin{equation}
\overline{\psi }_{\theta ,\delta }(x:\lambda ,\mu ),\underline{\psi }%
_{\theta ,\delta }(x:\lambda ,\mu )\rightarrow \overline{u}_{\lambda ,\mu
}(x),\hbox{ for any .}x\in [0,1]\hbox{, if }\theta \rightarrow 0\hbox{
and }\delta \rightarrow 0.  \label{continuity parameters}
\end{equation}
\end{proposition}

The following Figure explains the matching constructed for the auxiliary
barrier functions: \newline
\begin{center}
\includegraphics[scale=0.35]{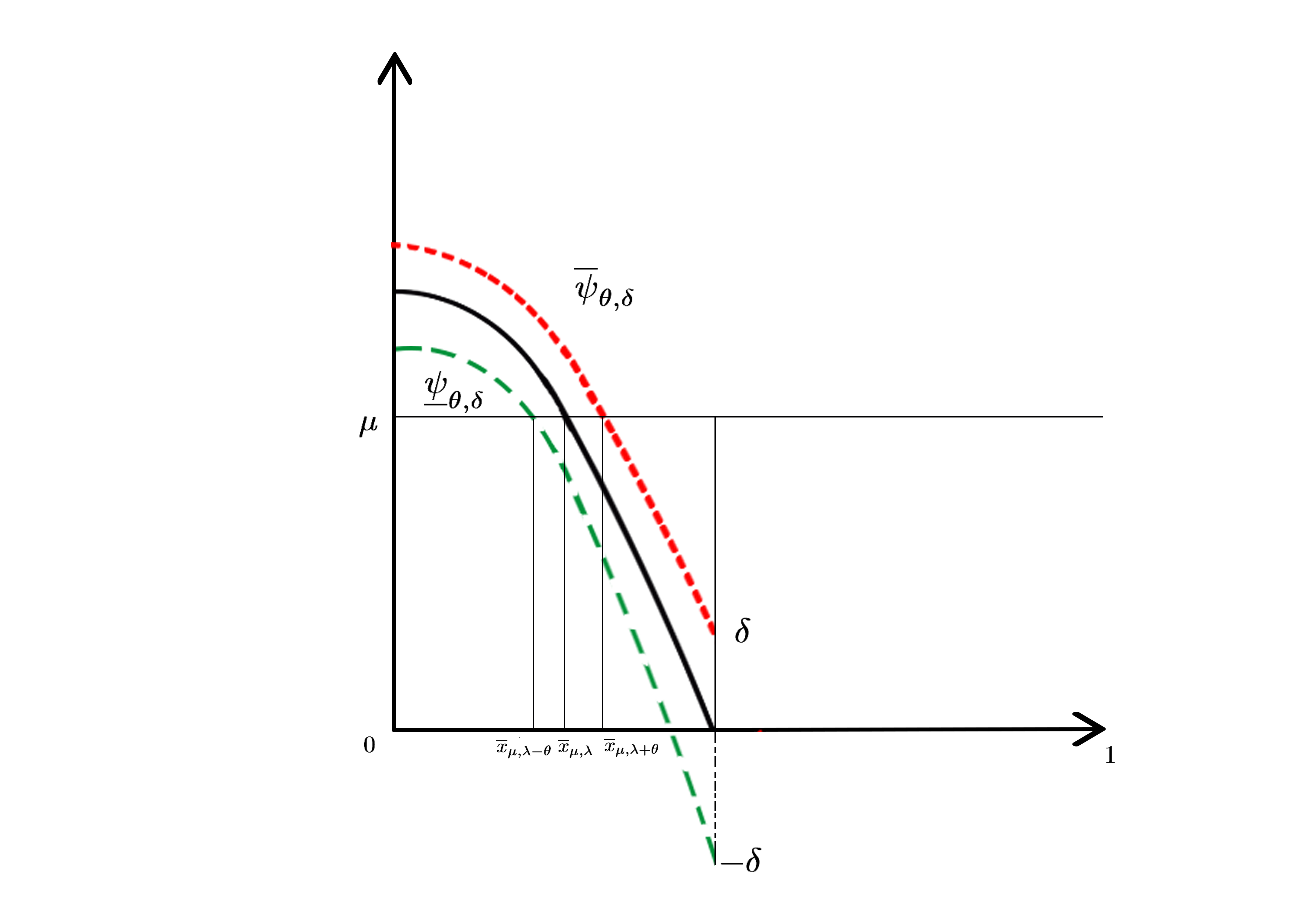}
\captionof{figure}{Construction of the auxiliary barrier functions. }
\end{center}
On the other hand, we shall prove that the set $E_{\theta ,\delta }:=\{v\in
L^{\infty }(0,1),$ $\underline{\psi }_{\theta ,\delta }<v<\overline{\psi }%
_{\theta ,\delta }\}$ is a $L^{\infty }-$attractive set with respect to the
dynamic problem if the parameters $\theta $ and $\delta $ are small enough:

\begin{proposition}
For some positive parameters $\theta $ and $\delta $, there exist $%
\varepsilon =\varepsilon (\theta ,\delta )>0$, such that if $\parallel u_{0}-%
\overline{u}_{\lambda ,\mu }\parallel _{L^{\infty }}<\varepsilon ,$ then $%
u(t,.:u_{0})\in $ $E_{\theta ,\delta }$ for any $t>0.$
\end{proposition}

\noindent \textit{Proof of Proposition 3.1. }For $\theta >0,$ small enough,
and let $\overline{x}_{\lambda +\theta ,\mu }$ the free boundary associated
to the parameter $\lambda +\theta $. \ By the continuity of the function (%
\ref{function g}) we know that there exists a small $\overline{h}(\theta )>0$
such that $\overline{x}_{\lambda +\theta ,\mu }=\overline{x}_{\lambda ,\mu }+%
\overline{h}(\theta )$. Let us\textit{\ }construct\textit{\ }the upper
barrier function\textit{\ }$\overline{\psi }_{\theta ,\delta }$ in the
following way: \ $\overline{\psi }_{\theta ,\delta }$ satisfies two
different boundary value problems over different regions.
\[
(\overline{P}_{L})\left\{
\begin{array}{ll}
-\overline{\psi }_{\theta ,\delta }^{\prime \prime }(x)+\omega ^{2}\overline{%
\psi }_{\theta ,\delta }(x)=\lambda +\theta & \quad \mbox{
in }\ (0,\overline{x}_{\lambda +\theta ,\mu }), \\[0.3cm]
\overline{\psi }_{\theta ,\delta }^{\prime }(0)=0,\hbox{ }\overline{\psi }%
_{\theta ,\delta }(\overline{x}_{\lambda +\theta ,\mu })=\mu . & \quad%
\end{array}%
\right.
\]
Moreover, for $\delta >0,$ small enough, we construct $\overline{\psi }%
_{\theta ,\delta }$ on $(\overline{x}_{\lambda +\theta ,\mu },1)$ such that
\[
(\overline{P}_{R})\left\{
\begin{array}{ll}
-\overline{\psi }_{\theta ,\delta }^{\prime \prime }(x)+\omega ^{2}\overline{%
\psi }_{\theta ,\delta }(x)=(\lambda +\theta )f_{0} & \quad \mbox{
in }\ (\overline{x}_{\lambda +\theta ,\mu },1), \\[0.3cm]
\overline{\psi }_{\theta ,\delta }(\overline{x}_{\lambda +\theta ,\mu })=\mu
,\hbox{ }\hbox{ }\overline{\psi }_{\theta ,\delta }(1)=\delta . & \quad%
\end{array}%
\right.
\]%
Thus, from the strong maximum principle for linear equations and from the
above construction of problems $(\overline{P}_{L})$ and $(\overline{P}_{R}),$
we have
\[
\overline{\psi }_{\theta ,\delta }(x)>\overline{u}_{\lambda ,\mu }(x)\quad %
\hbox{in}\quad (0,1).
\]%
Moreover, it is clear that from the study of $\overline{u}_{\lambda ,\mu
}(x) $ made in the above Section that $\overline{\psi }_{\theta ,\delta
}(x:\lambda ,\mu )\rightarrow \overline{u}_{\lambda ,\mu }(x),$ for any $%
x\in (0,1)$ if $\theta \rightarrow 0$ and $\delta \rightarrow 0.$ Notice
also that from the uniform continuity of $\overline{\psi }_{\theta ,\delta }$
on $[0,1]$ there exists $\varepsilon =\varepsilon (\theta ,\delta )>0$,
small enough,$\ $such that
\begin{equation}
0<\overline{\psi }_{\theta ,\delta }(x)-\overline{u}_{\lambda ,\mu }(x)\leq
\varepsilon \hbox{ for any }x\in \lbrack 0,1].  \label{difference epsilon +}
\end{equation}%
Notice that the easier choice $\overline{\psi }_{\theta ,\delta }(x)=$\ $%
\overline{u}_{\lambda +\theta ,\mu }(x)$ is not good (for our purposes),
near $x=1$, since $\overline{u}_{\lambda +\theta ,\mu }(1)=\overline{u}%
_{\lambda ,\mu }(1)=0$ and we shall need later a strict positive distance in
$L^{\infty }(0,1).$

\noindent In a similar way, the continuity of the function (\ref{function g}%
) implies that there exists a small $\underline{h}(\theta )>0$ such that \ $%
\overline{x}_{\lambda -\theta ,\mu }=\overline{x}_{\lambda ,\mu }-\underline{%
h}(\theta ).$ The lower barrier function $\underline{\psi }_{\theta ,\delta
} $ can be built, then, by means of the following auxiliary problems:
\[
(\underline{P}_{L})\left\{
\begin{array}{ll}
-\underline{\psi }_{\theta ,\delta }^{\prime \prime }(x)+\omega ^{2}%
\underline{\psi }_{\theta ,\delta }(x)=\lambda -\theta & \quad \mbox{
in }\ (0,\overline{x}_{\lambda -\theta ,\mu }), \\[0.3cm]
\underline{\psi }_{\theta ,\delta }^{\prime }(0)=0,\hbox{ }\underline{\psi }%
_{\theta ,\delta }(\overline{x}_{\lambda -\theta ,\mu })=\mu , & \quad%
\end{array}%
\right.
\]%
\[
(\underline{P}_{R})\left\{
\begin{array}{ll}
-\underline{\psi }_{\theta ,\delta }^{\prime \prime }(x)+\omega ^{2}%
\underline{\psi }_{\theta ,\delta }(x)=(\lambda -\theta )f_{0} & \quad
\mbox{
in }\ (\overline{x}_{\lambda -\theta ,\mu },1), \\[0.3cm]
\underline{\psi }_{\theta ,\delta }(\overline{x}_{\lambda -\theta ,\mu
})=\mu ,\hbox{ }\hbox{ }\underline{\psi }_{\theta ,\delta }(1)=-\delta . &
\quad%
\end{array}%
\right.
\]%
Using again the strong maximum principle for linear equations, we get
\[
\underline{\psi }_{\theta ,\delta }(x)<\overline{u}_{\lambda ,\mu }(x)\quad %
\hbox{in}\quad (0,1).
\]%
As before, from the study of $\overline{u}_{\lambda ,\mu }(x)$ made in the
above Section, $\underline{\psi }_{\theta ,\delta }(x:\lambda ,\mu
)\rightarrow \overline{u}_{\lambda ,\mu }(x),$ for any $x\in (0,1)$ if $%
\theta \rightarrow 0$ and $\delta \rightarrow 0.$ We point out that, \ as a
matter of fact, $0\leq \underline{\psi }_{\theta ,\delta }(x)\leq -\delta $
for any $x\in \lbrack 1-\rho (\theta ,\delta ),1],$ for some $\rho (\theta
,\delta )>0$ with $\rho (\theta ,\delta )\rightarrow 0$ if $\theta
\rightarrow 0$ and $\delta \rightarrow 0.$ Moreover, from the uniform
continuity of $\underline{\psi }_{\theta ,\delta }$ on $[0,1]$ there exists $%
\varepsilon =\varepsilon (\theta ,\delta )>0$, small enough,$\ $such that
\begin{equation}
-\varepsilon <\overline{\psi }_{\theta ,\delta }(x)-\overline{u}_{\lambda
,\mu }(x)<0\hbox{ for any }x\in \lbrack 0,1]._{\blacksquare }
\label{difference epsilon -}
\end{equation}

\noindent \textit{Proof of Proposition 3.2 \ }Even if $\omega ^{2}>0,$ the
proof of this property is entirely similar to the correspondent part of the
proof of Theorem 1.3 of \cite{BensidDiaz} (concerning the case $\omega =0$).
We send the reader to this paper for the details$._{\blacksquare }$\newline
\newline
\noindent \textit{Proof of Theorem 1.2 } To conclude that $\overline{u}%
_{\lambda ,\mu }$ is stable in $L^{\infty }(0,1),$ it suffice to combine
Propositions 3.1 and 3.2 and to use that the invariant set $E_{\theta
,\delta }$ contains the $L^{\infty }(0,1)-$neighborhood of $\overline{u}%
_{\lambda ,\mu }$ of radium $\varepsilon =\varepsilon (\theta ,\delta ).$ $%
_{\blacksquare }.$

\section{Instability of the decreasing part of the bifurcation curve}

In this section, we shall prove Theorem 1.3 and Theorem 1.4 concerning the
instability of lower branch $\underline{u}_{\lambda ,\mu }.$\newline

\bigskip

\noindent \textit{Proof of Theorem 1.3 \ }It is based on the fact that if
the principal eigenvalue of $P_{\eta }(x_{\mu ,\lambda }:f_{0},\lambda )$ is
negative then the stationary monotone solution is unstable. The main idea of
the proof is an adaptation of the similar result presented in \cite%
{BensidDiaz} for the case $\omega =0$. If we consider the solutions $v$ of
the parabolic problem $PP^{\ast }(\lambda ,\beta ,v_{0})$\ with $%
v_{0}=u_{\lambda }+\phi w_{0}$ with $w_{0}$ smooth and $\phi $ small enough
\ and to approximate $v$, as $\phi \rightarrow 0$, by functions of the form $%
v(x,t)=u_{\lambda }(x)+\phi w(x,t)$ with $w(t,x)=e^{-\nu t}U(x)$, with $U$
solution of the eigenvalue problem $P_{\eta }(x_{\lambda ,\mu
}:f_{0},\lambda ).$ Thus, using Proposition 4.1 in \cite{BensidDiaz}, we
find that $U$ satisfies
\[
\left\{
\begin{array}{l}
-U^{\prime \prime }=\rho U,%
\hbox{ \ \ \ \ \ \ \ \ \ \ \ \ \ \ \ \ \ \ \ \ \
\ \ \ \ \ \ \ \ }x\in (0,x_{\lambda ,\mu })\cup (x_{\lambda ,\mu },1), \\%
[0.3cm]
U^{\prime }(0)=0,\quad U(1)=0, \\[0.3cm]
U_{-}(x_{\lambda ,\mu })=U_{+}(x_{\lambda ,\mu }),\hbox{ }U_{-}^{\prime
}(x_{\lambda ,\mu })-U_{+}^{\prime }(x_{\lambda ,\mu })=\lambda
(1-f_{0})U(x_{\lambda ,\mu }),%
\end{array}%
\right.
\]%
where $\rho :=\eta -\omega ^{2}.$ We recall that in \cite{BensidDiaz}, we
have studied the case $\omega =0.$ In this case and when $\rho =\eta :=-\tau
^{2},$ we have showed that the free boundary $\underline{x}_{\lambda ,\mu
}\in (0,\frac{1}{2-f_{0}})$ generate $\tau >0$ which gives a positive
solution $U.$\newline
When $\omega \neq 0$ but $0<\omega <1,$ the same techniques can be used here
for $\rho =\eta -w^{2}:=\tau ^{\ast }$ to prove the existence of $\tau
^{\ast }>0$ given by
\[
\tau ^{\ast }=\tau ^{2}+\omega ^{2}.
\]%
This choice is always possible since $\omega ^{2}\in (0,1)._{\blacksquare }$
\newline

\bigskip

\noindent \textit{Proof of Theorem 1.4. }We shall use a very special
construction of auxiliary initial data which leads, after a suitable
convergence, to solution $\underline{u}_{\lambda ,\mu }$. Let $\theta >0$
small enough and take
\[
\overline{u}_{0}(x):=\underline{u}_{\lambda -\theta ,\mu }(x).
\]%
From the the proof of the results iii) and v) of Theorem 1.1, we have that $%
\theta >0$ implies
\[
\overline{u}_{0}(x)>\underline{u}_{\lambda ,\mu }(x),\quad \forall x\in
\lbrack 0,1).
\]%
Moreover, to check condition (\ref{epsilon near}) observe that it is clear
that $\left\Vert \underline{u}_{\lambda -\theta ,\mu }\right\Vert
_{L^{\infty }(0,1)}$ depends continuously with respect to $\theta $. Indeed,
this is exactly the continuity condition of the bifurcation curve given in
Theorem 1.1. Notice, for instance, that the maximum point of $\underline{u}%
_{\lambda -\theta ,\mu }$ takes place at $x=0$ and that, although at this
point we only have a Neumann boundary condition, the fact that $\underline{u}%
_{\lambda -\theta ,\mu }$ verifies
\begin{equation}
\left\{
\begin{array}{ll}
-\underline{u}_{\lambda -\theta ,\mu }^{\prime \prime }+\omega ^{2}%
\underline{u}_{\lambda -\theta ,\mu }(x)=\lambda -\theta , & \quad
\mbox{
for}\ x\in (0,\underline{x}_{\lambda -\theta ,\mu }) \\[0.3cm]
\underline{u}_{\lambda -\theta ,\mu }^{\prime }(0)=0,\hbox{ },\underline{u}%
_{\lambda -\theta ,\mu }(\underline{x}_{\lambda -\theta ,\mu })=\mu , & \quad%
\end{array}%
\right.  \label{problem theta}
\end{equation}%
implies the above mentioned continuous dependence as a by-product of the
continuous dependence of solutions of problem (\ref{problem theta}) with
respect to the $L^{\infty }(0,1)$ norm of the right hand function and the
continuous dependence with respect to the own interval of definition (recall
that we know that the continuity of the function (\ref{function g}) implies
that there exists a continuous function $\underline{h}(\theta )>0$ such that
\ $\underline{x}_{\lambda -\theta ,\mu }=\underline{x}_{\lambda ,\mu }-%
\underline{h}(\theta )$). \ Thus we have that, given $\theta >0$ small
enough, there exists $\varepsilon =\varepsilon (\theta )>0$ such that
\begin{equation}
0<\overline{u}_{0}(x)-\underline{u}_{\lambda ,\mu }(x)\leq \varepsilon \quad
\forall x\in \lbrack 0,1),  \label{entorno en L infinito}
\end{equation}%
which shows (\ref{epsilon near}).\newline
\newline
\noindent To complete the proof of Theorem 1.4 we also need the following
auxiliary result:

\begin{lemma}
If $u_{0}\in W^{2,\infty }(0,1)$ is non-degenerate and
\[
(P_{0})\left\{
\begin{array}{ll}
-u_{0}^{\prime \prime }+\omega ^{2}u_{0}(x)\leq \lambda \beta (u_{0}) &
\quad \mbox{
a.e }\ x\in (0,1), \\[0.3cm]
u_{0}(0)=0,\hbox{ }u_{0}^{\prime }(0)\leq 0, & \quad%
\end{array}%
\right.
\]%
then the unique non-degenerate solution of $PP^{\ast }(\lambda ,\beta
,u_{0}) $ satisfies that
\begin{equation}
\frac{\partial u}{\partial t}(t,x)\geq 0\quad \hbox{a.e}\quad x\in
(0,1)\quad \hbox{and a.e}\quad t>0.
\end{equation}
\end{lemma}

\noindent \textit{Proof of Lemma 4.1. } Let $\beta _{\varepsilon }(.)$ a
sequence of smooth increasing functions such that $\beta _{\varepsilon
}(.)\rightarrow \beta $ in the sense of $\mathbb{R}^{2}$ maximal monotone
graphs. Then, if $u_{\varepsilon }$ is the unique solution of $PP^{\ast
}(\lambda ,\beta _{\varepsilon },u_{0}),$ we have that
\[
u_{\varepsilon }\rightarrow u\quad \hbox{in}\quad C([0,T]:L^{1}(0,1))\quad %
\hbox{as}\quad \varepsilon \rightarrow 0,
\]%
with $u$ the unique non-degenerate solution of $PP^{\ast }(\lambda ,\beta
,u_{0})$ (see, e.g. \cite{Gianni})$.$ Then, it is enough to prove $(4.16)$
for $u_{\varepsilon }.$ Indeed, since
\begin{equation}
\frac{\partial u_{\varepsilon }}{\partial t}(t,x)\geq 0\Leftrightarrow
u_{\varepsilon }(t,.)\leq u_{\varepsilon }(t^{\prime },.),\quad \quad
\forall t<t^{\prime },
\end{equation}%
this implies (as $\varepsilon \rightarrow 0$)
\[
u(t,.)\leq u(t^{\prime },.),\quad \forall t<t^{\prime },
\]%
and thus we have $(4.16).$\newline
Finally, taking $v_{\varepsilon }:=\frac{\partial u_{\varepsilon }}{\partial
t}.$ By differentiation in the PDE, the function $v_{\varepsilon }$ verifies
\[
DPP^{\ast }(\lambda ,\beta ,v_{0})\left\{
\begin{array}{ll}
v_{t}-v_{xx}+\omega ^{2}v=\lambda \beta _{\varepsilon }^{\prime
}(u_{\varepsilon })v & \quad x\in (0,1),\hbox{ }t>0, \\[0.3cm]
v_{x}(0,t)=0,\quad v(1,t)=0 & \quad t>0, \\[0.3cm]
v(x,0)=\frac{\partial u}{\partial t}(x,t)|_{t=0}, & \quad x\in (0,1),%
\end{array}%
\right.
\]%
but
\[
\frac{\partial u}{\partial t}(.,t)|_{t=0}=u(.,t)_{xx}-\omega
^{2}u(.,t)+\lambda \beta _{\varepsilon }(u(.,t))|_{t=0}\geq 0,\quad
\]%
thanks to the assumption $(P_{0}).$ Then, by the maximum principle for $%
DPP^{\ast }(\lambda ,\beta ,v_{0}),$ we conclude that $v\geq 0,$ for $%
(t,x)\in \lbrack 0,+\infty )\times (0,1)$ and so we have $%
(4.17)._{\blacksquare }$\newline
\newline
\textit{An alternative proof of Lemma 4.1}. Let $v=\frac{\partial u}{%
\partial t},$ then $v\in C^{0}((0,+\infty )\times (0,1))$ and verifies
\begin{equation}
\left\{
\begin{array}{ll}
v_{t}-v_{xx}+\omega ^{2}v=\lambda \delta _{\{x(t)\}}v & \quad x\in (0,1),%
\hbox{ }t>0, \\[0.3cm]
v_{x}(0,t)=0,\quad v(1,t)=0 & \quad t>0, \\[0.3cm]
v(x,0)=\frac{\partial u}{\partial t}(x,t)|_{t=0}\geq 0. & \quad x\in (0,1),%
\end{array}%
\right.
\end{equation}%
The maximum principle applies to $(4.18)$ and we can conclude that $v\geq 0$
(problems with measures of similar nature in \cite{DiazFowler} or \cite%
{Liang-Lin-Matano}). $_{\blacksquare }$\newline
\newline
Now, to conclude the proof of Theorem 1.4, we remark that $(1.4)$ is already
shown in Lemma 4.1. To see $(1.5),$ it is enough to check that $\overline{u}%
_{\lambda ,\mu }$ verifies $PP^{\ast }(\lambda ,\beta ,\overline{u}%
_{\varepsilon })$ and that
\[
\underline{u}_{\lambda -\theta ,\mu }(x)=\overline{u}_{0}(x)\leq \overline{u}%
_{\lambda ,\mu }(x)\hbox{ on }(0,1).
\]%
Thus, since $u(t:u_{0})$ and $u_{\lambda ,\mu }(.)$ are non degenerate (see
\cite{BensidDiaz}) and using the comparison principle for non-denenerate
solutions of $PP^{\ast }(\lambda ,\beta ,v_{0}),$ we have that
\[
u(t:u_{0})\leq \overline{u}_{\lambda ,\mu },\quad \quad \forall t>0%
\hbox{ on
}(0,1).
\]%
Finally, if we define the sequence $\{\Vert u(t_{n})\Vert
_{L^{p}(0,1)}\}_{t_{n}>0}$, $\forall p>1,$ then by Lemma 4.2, we see that
this sequence is increasing and bounded. Hence, using a classical result we
conclude that there exists $\xi \in \mathbb{R},$ with $\xi \leq \Vert
\overline{u}_{\lambda ,\mu }\Vert _{L^{p}(0,1)}$ such that $%
u(t_{n})\rightharpoonup u_{\infty }$ where $u_{\infty }$ is a weak solution
of the stationary problem corresponding to
\[
-u_{\infty }^{\prime \prime }+\omega ^{2}u_{\infty }=\lambda f(u_{\infty })%
\hbox{ on }(0,1).
\]%
\newline
In fact, by Theorem 1 of \cite{Diaz-Her-Tello}, we have $(1.6).$\newline
A similar conclusion (proving (\ref{convergence u star})) can be obtained
for the solution $u(t:\underline{u}_{0})$ once we choose $\underline{u}%
_{0}(x):=\underline{u}_{\lambda +\theta ,\mu }(x).$ $_{\blacksquare }$

\bigskip

The following figure explains the dynamics of some initial data which are in
small $L^{\infty }(0,1)-$neighborhood of the instable solution $\underline{u}%
_{\lambda ,\mu }:$

\begin{center}
\includegraphics[scale=0.3]{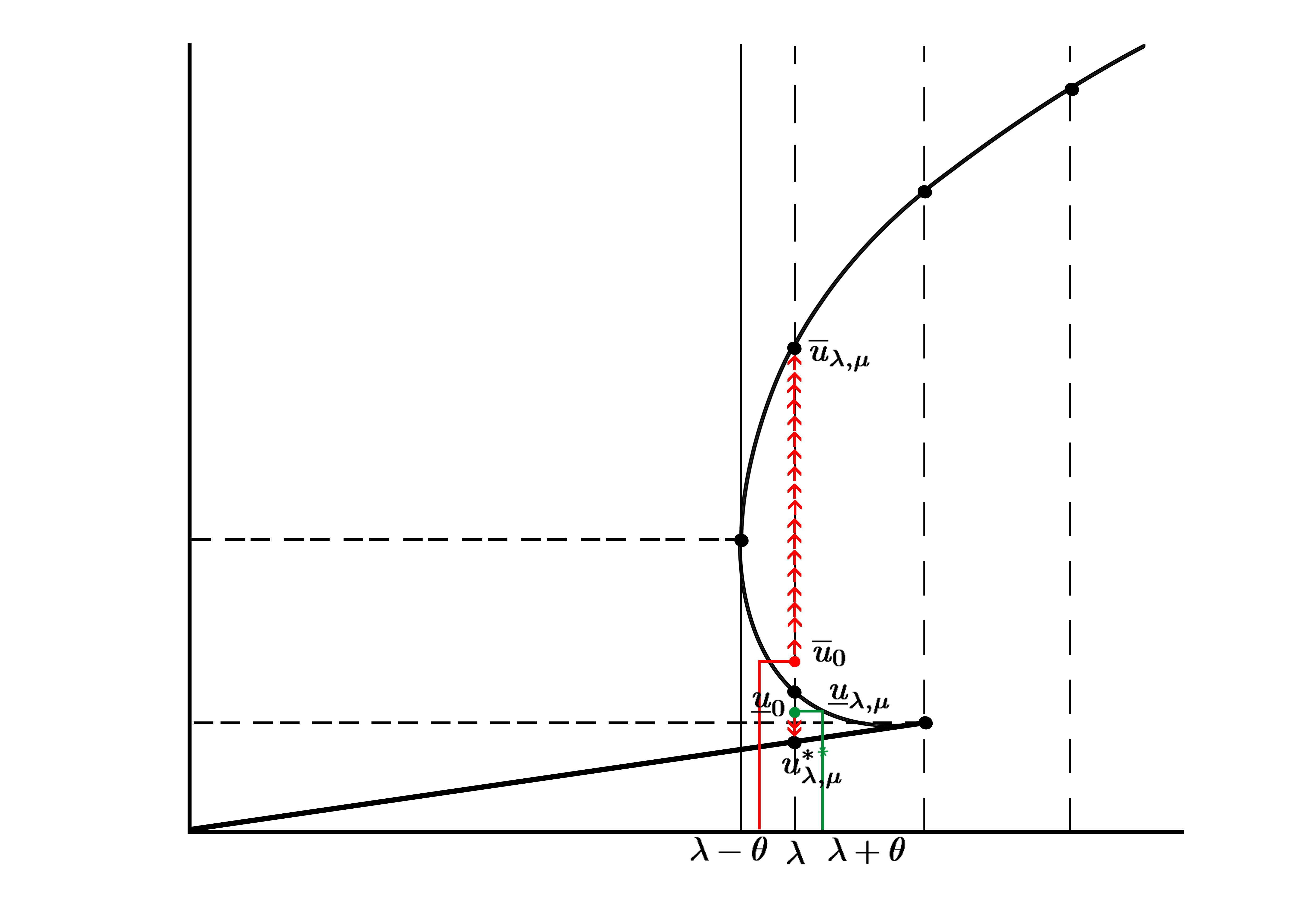}
\captionof{figure}{Dynamics of solutions corresponding to suitable initial data
closed to the unstable equilibrium $\underline{u}_{\lambda ,\mu }.$ }
\end{center}

\begin{remark}
The proof of Theorem 1.4 also shows that the transient free boundaries
corresponding to those initial data satisfy that $\underline{x}_{\mu
,\lambda }(t)\nearrow \overline{x}_{\lambda ,\mu }<1$ as $t\rightarrow
+\infty $ and that $\underline{x}_{\lambda ,\mu }(t)\searrow 0$ as $%
t\rightarrow +\infty .$\newline
\end{remark}

\begin{remark}
The above proof generalizes and improves (with a different point of view)
the results of \cite{Issard} for the case of $\beta $ a regular function.
\end{remark}

\bigskip

\section*{Acknowledgments}

We thank to one anonymous referee for the received constructive comments
after a very careful reading of the previous version of this manuscript.
This research has partially been supported by the MINECO (Spain) project
MTM2014-57113-P of the DGISPI (Spain) and the Research Group MOMAT (Ref.
910480) supported by the Universidad Complutense de Madrid.

\begin{center}
{\footnotesize
\begin{tabular}{ll}
Sabri Bensid & Jes\'{u}s Ildefonso D\'{\i}az \\
Dynamical Systems and Applications Laboratory & Instituto de Matem\'{a}tica
Interdisciplinar \\
Department of Mathematics, Faculty of Sciences & Depto. de An\'{a}lisis Matem%
\'{a}tico-Matem\'{a}tica Aplicada \\
University of Tlemcen, B.P. 119 & Universidad Complutense de Madrid \\
Tlemcen 13000 & Plaza de las Ciencias 3, 28040--Madrid \\
Algeria & Spain \\
edp\_sabri@yahoo.fr & jidiaz@ucm.es%
\end{tabular}%
}
\end{center}

\end{document}